\documentclass[12pt]{amsart}
\usepackage{amssymb,amsmath,amstext,amsthm,graphicx,setspace,color}
\usepackage{amsmath,amssymb,amscd,amsmath,amstext,mathrsfs,amsfonts,pb-diagram,algorithmic}%

\newenvironment{Macaulay2}{
\begin{spacing}{0.7}
\smallskip
\small
} {
\smallskip
\end{spacing}
\medskip
}

\theoremstyle{plain}

\newtheorem{algorithm}{Algorithm}

\newtheorem{conjecture}{Conjecture}

\newtheorem{definition}{Definition}

\numberwithin{equation}{section}

\def\NAG{{\em NAG4M2}}
\def\M2{{\em Macaulay2}}
\def\Mathematica{{\em Mathematica}}

\newcommand{\CHd}{\mathcal{H}_{(d)}}

\newcommand{\Pn}{\P(\C^{n+1})}

\newcommand{\Real}{\operatorname{Re}}

\def\C{{\mathbb C}}

\def\Q{{\mathbb Q}}

\def\P{{\mathbb P}}
\def\S{{\mathbb S}}

\def \LinearHomotopyConstant{0.04804448}

\let\oldmarginpar\marginpar
\renewcommand\marginpar[1]{\-\oldmarginpar[\raggedleft\footnotesize #1]%
{\raggedright\footnotesize #1}}

\title[A search for an optimal start system]{A search for an optimal start system\\for numerical homotopy continuation}
\author{Anton Leykin}
        \thanks{Anton Leykin. School of Mathematics, Georgia Tech, Atlanta GA, USA ({\tt leykin@math.gatech.edu}).
         Partially supported by NSF grant DMS-0914802}

\date{\today}

\begin{document}
\begin{abstract}
We use our recent implementation of a certified homotopy tracking algorithm to search for start systems that minimize the average complexity of finding all roots of a regular system of polynomial equations. While finding optimal start systems is a hard problem, our experiments show that it is possible to find start systems that deliver better average complexity than the ones that are commonly used in the existing homotopy continuation software.
\end{abstract}
\maketitle

\section{Introduction}
This section introduces homotopy methods for solving polynomial systems and sets up the stage for the experiments that are in the core of this article. The main contribution of this work concerns optimal start systems, which do not enter this stage until Subsection~\ref{subsec: intro to optimal systems}.

\smallskip

The problem of finding a good algorithm to solve a system of polynomial equations has been on the minds of many mathematicians for many centuries. Here we address this problem in its most basic form:
\begin{itemize}
  \item the coefficients of the polynomials are complex numbers;
  \item the system is square (neither under- nor overdetermined);
  \item the system has finitely many complex solutions and this number is equal to the the number of solutions of a generic system of equations with the same degrees (the B\'ezout bound) and, therefore, all solutions are regular.
\end{itemize}

In principle, the stated problem can be solved symbolically. Let field $C$ be an extension of $\Q$ containing all coefficients. Then, for instance, via elimination techniques employing Gr\"obner bases one can find a finite extension of $C$ that contains all coordinates of all solutions.
In practice, for systems with a large number of solutions this method could be impractical: the complexity of expressions representing the resulting extension and solutions is high, not to mention the phenomena of intermediate expression swell. Besides, even if these expressions are obtained, getting numerical approximations to the solutions could be a nontrivial task.

Approximate solutions would arguably constitute the most widely applicable output of any polynomial system solving algorithm. A basic question to ask at this point is: What does it mean for a solution to be approximate? Here are two answers:
\begin{itemize}
  \item fix $\varepsilon>0$, then $z\in\C^n$ is an approximation of an exact solution $\zeta \in \C^n$ if $\|z-\zeta\|\leq \varepsilon$;
  \item fix a refinement algorithm, i.e.: a function $R:\C^n\to\C^n$, then $z$ is an approximate zero associated to $\zeta$ if the sequence $\{ R^k(z) \}$ converges to $\zeta$.
\end{itemize}
The first way, while frequently preferred by practitioners, is clearly dependent on the threshold $\varepsilon$. On top of that, one relaxation is commonly made in applications: the distance $\|z-\zeta\|$ is replaced by its estimate.

The second approach is intrinsic once the refinement procedure $R$ is fixed. Set $R$ to be Newton's iterator to get the notion of Smale's {\em approximate zero} \cite{Smale:approximate-zero} used as a basis for $\alpha$-theory. (These concepts are defined in Section~\ref{sec: preliminaries}.) It is this meaning of an approximate solution that we adopt for the purpose of this article.

For the purpose of finding approximate solutions let us introduce the homotopy continuation technique: given a {\em target system} $f$ construct a {\em start system} $g$ and track the curve
$$t\rightarrow h_t,\ t\in[0,T]; \text{ such that } h_T=f \text{ and } h_0=g;$$
in the space of polynomial systems. An atomic task of a homotopy continuation algorithm is to take a {\em start solution}, i.e., a zero $z_0 \in g^{-1}(0)$, and track the homotopy path starting at $z_0$ and finishing at a {\em target solution}, a zero of $f$. For details on how $g$ and $t\rightarrow h_t$ are set up so that the last sentence makes sense one may consult the recent book~\cite{Sommese-Wampler-book-05}; we construct several particular homotopies in Section~\ref{sec: preliminaries}.

In the last decade several software systems (e.g.: \cite{Bertini}, \cite{HOM4PSwww}, \cite{Leykin:NAG4M2}, \cite{V99}) based on homotopy continuation powered by numerical predictor-corrector methods have been developed. Problems with millions of solutions could be attacked due to speed achieved through path-tracking heuristics. The target solutions in some cases can be {\em certified} rigorously, however, usually only data characterizing the {\em quality} of a solution (such as estimates for condition number, error, residual value, etc.) is produced. The latter is aimed at providing the software user with {\em confidence} in the obtained results, which is arguably what the most users need.

\subsection{Certified numerical homotopy tracking} 
\begin{figure}[h]
\begin{picture}(300,220)
   \put(10,10){\includegraphics[scale=0.85, keepaspectratio]{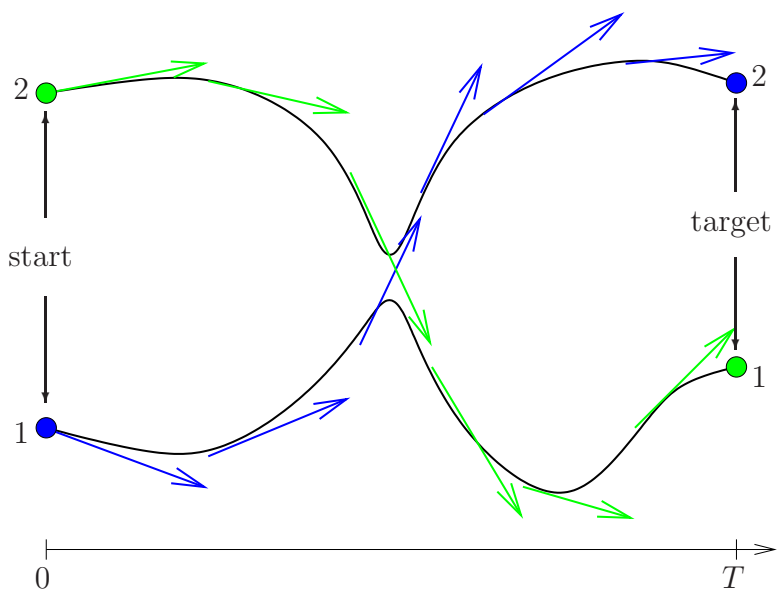}}
   \put(10,0){$0$}
   \put(2,55){$1$}
   \put(2,185){$2$}
   \put(0,122){start}
   \put(14,140){\vector(0,1){40}}
   \put(14,110){\vector(0,-1){40}}
   \put(270,0){$T$}
   \put(281,76){$1$}
   \put(281,190){$2$}
   \put(258,135){target}
   \put(275,150){\vector(0,1){35}}
   \put(275,125){\vector(0,-1){35}}
\end{picture}
  \caption{{``Path-crossing'' scenario: the tracker jumps from path 1 to path 2 and vice-versa.}}
  \label{fig: path switching}
\end{figure}

Apart from the certification of the end points of the homotopy paths, we would also like to certify that {\em all points} produced by the tracker working on one homotopy path are approximate zeros associated to the points on this path. Note that a heuristic tracker could lead to the picture in Figure~\ref{fig: path switching}: at every point where a prediction is evaluated (depicted by arrows) the system is well-conditioned, approximate solutions to all solutions of the target system are discovered, however the continuation paths are swapped by the tracker.

Ensuring this does not happen is of ultimate importance to the algorithms that produce discrete output based on the assumption that no {\em path-jumping} occurs. For instance, discovering monodromy via homotopy tracking is a part of algorithms for the genus of an irreducible curve within an algebraic set \cite{Bates-Peterson-Sommese-Wampler:genus}, Galois groups of Schubert problems \cite{Leykin-Sottile:HoG}, and numerical irreducible decomposition. All these are applications to problems in pure mathematics and certification is highly desirable in order to give the obtained results the status of a mathematical {\em proof}.

Creation of a certified homotopy tracking procedure was the goal of our joint work with Carlos Beltr\'an in \cite{Beltran-Leykin:CHT}: the resulting algorithm has been implemented in the {\em NumericalAlgebraicGeometry} package \cite{Leykin:NAG4M2} of \M2\ \cite{M2www}, referred to as \NAG.

\subsection{Start systems with optimal average complexity} \label{subsec: intro to optimal systems}
We used the implementation mentioned above to estimate the average complexity of computation with the linear homotopy using various {\em initial pairs} (start system, start solution). The results concerning the average complexity of finding {\em one} solution of a system obtained in~\cite{Beltran-Leykin:CHT} are outlined in subsection~\ref{subsec: previous experiments}.

The main question considered in this paper and not addressed in~\cite{Beltran-Leykin:CHT} is: How to find an {\em optimal} start system, i.e., one that delivers the minimal average complexity of computing {\em all} solutions?\footnote{In the context of structured polynomial systems, the word {\em optimal} is used in a different sense: it refers to a homotopy that minimizes the number of continuation paths tracked in the computation. For example, polyhedral homotopies implemented in~\cite{V99, HOM4PSwww} are optimal for solving sparse polynomial systems in that sense. The exact meaning of {\em optimality} in the context of this paper shall become clear in Section~\ref{sec: preliminaries}.}

This is a hard question, which has not been answered even in the one-equation case. The work of Shub and Smale~\cite{ShubSmale:Bezout-3} suggests that an optimal start system for the linear homotopy has to have the smallest {\em condition number}. They define the concept of the condition number for a polynomial equation and prove that there is a family of polynomials --- elliptic Fekete polynomials --- for which the condition number grows polynomially with the degree.

Subsection~\ref{subsec: previous experiments} revisits the experiments of~\cite{Beltran-Leykin:CHT} with a view towards a search of an optimal initial pair (for finding {\em one} solution of a target system). In subsection \ref{subsec: new experiments} the results of several new experiments are presented in detail show that it is possible to find start systems that give better average complexity (of finding {\em all} roots) than the total-degree homotopy systems that are a popular choice in current software.

While the discovery and experiments are done with the certified homotopy tracking algorithm, heuristic computation can benefit from finding either optimal or near-optimal start systems as well. This is backed up by the data obtained by running an implementation of a heuristic algorithm in several experiments.

The \M2~\cite{M2www} and \Mathematica~\cite{MATHEMATICAwww} scripts used in this paper are published at {\tt http://people.math.gatech.edu/\~{}aleykin3/OSS/}.

\subsection{Acknowledgements} The author would like to thank the referees, the organizers of the BIRS workshop on ``Randomization, Relaxation, and Complexity'' in 2010, and Institut Mittag-Leffler where the work on this article has been completed.

\section{Preliminaries}\label{sec: preliminaries}

To simplify the discussion in this paper we shall always consider the homogeneous problem of polynomial system solving. For a positive integer $l\geq1$, let $\mathcal{H}_l$ be the vector space of all homogeneous polynomials of degree $l$ with complex coefficients and unknowns $X_0,\ldots,X_n$. For a list of degrees $(d)=(d_1,\ldots,d_n)$ let $$\CHd=\mathcal{H}_{d_1}\times\cdots\times\mathcal{H}_{d_n}.$$
An element $h\in\CHd$ is seen as a vector in a space of dimension $$N + 1 = \sum_{i=1}^n\binom{n+d_i}{d_i}$$
and as a system of $n$ homogeneous equations with $n+1$ unknowns. The zero set of $h$ is a projective variety: we consider its zeros as projective points $\zeta\in\Pn$. Throughout the paper no distinction between a point in $\Pn$ and a representative of the point in $\C^{n+1}$ is made: when necessary it is implied that a representative has the unit norm.

\subsection{Approximate zeros and the projective Newton's operator}

It has been pointed out in the introduction that describing the zeros of $h\in\CHd$ exactly is a hard task: one may ask for points which are $\varepsilon$--close to some zero, however, this is not an intrinsic concept. What we employ here is the concept of an approximate zero of \cite{Smale:approximate-zero}.

First we define the projective Newton's method as follows. Let $h\in\CHd$ and $z\in\Pn$. Then,
\[
N_\P(h)(z)=z-\left(Dh(z)\mid_{z^\perp}\right)^{-1}h(z),
\]
where $Dh(z)$ is the $n\times (n+1)$ Jacobian matrix of $h$ at $z\in\Pn$, and
\[
Dh(z)\mid_{z^\perp}
\]
is the restriction of the linear operator defined by $Dh(z):\C^{n+1}\rightarrow\C^n$ to the orthogonal complement $z^\perp$ to the one-dimensional space spanned by $z$. Hence, $N_\P(h)(z)$ is defined as long as $\left(Dh(z)\mid_{z^\perp}\right)^{-1}$, a linear operator from $\C^n$ to $z^\perp\subset \C^{n+1}$ of dimension $n$, is invertible.

Let $d_R$ be the Riemann distance in $\Pn$:
\[
d_R(z,z')=\arccos\frac{|\langle z,z'\rangle|}{\|z\|\,\|z'\|}\in[0,\pi/2],
\]
where $\langle \cdot,\cdot\rangle$ and $\|\cdot\|$ are the usual Hermitian product and norm in $\C^{n+1}$. Note that $d_R$ is well defined on $\Pn\times\Pn$.
\begin{definition}\label{def:projzppzero}
We say that $z\in\Pn$ is an approximate zero of $h\in\CHd$ with associated zero $\zeta\in\Pn$ if $N_\P(h)^l(z)$ is defined for all $l\geq0$ and
\[
d_R(N_\P(h)^l(z),\zeta)\leq\frac{d_R(z,\zeta)}{2^{2^{l}-1}},\;\;\;\;l\geq0.
\]
\end{definition}

\subsection{Bombieri-Weyl norm}\label{sec:BW}
Given two polynomials $v,w\in\mathcal{H}_l$,
\[
v=\sum_{\alpha_0+\ldots+\alpha_n=l}a_{\alpha_0,\ldots,\alpha_n}X_0^{\alpha_0}\cdots X_n^{\alpha_n},
\]
\[
w=\sum_{\alpha_0+\ldots+\alpha_n=l}b_{\alpha_0,\ldots,\alpha_n}X_0^{\alpha_0}\cdots X_n^{\alpha_n},
\]
we define their Bombieri-Weyl product to be
\[
\langle v,w\rangle=\sum_{\alpha_0+\alpha_1+\ldots+\alpha_n=l}\binom{l}{(\alpha_0,\ldots,\alpha_n)}^{-1} a_{\alpha_0,\ldots,\alpha_n}\overline{b_{\alpha_0,\ldots,\alpha_n}},
\]
where $\overline{\;\cdot\;}\,$ is the complex conjugation and
\[
\binom{l}{(\alpha_0,\ldots,\alpha_n)}=\frac{l!}{\alpha_0!\cdots\alpha_n!}
\]
is the multinomial coefficient.

Then, given two elements $h=(h_1,\ldots,h_n)$ and $h'=(h'_1,\ldots,h'_n)$ of $\CHd$, we define
\[
\langle h,h'\rangle=\langle h_1,h_1'\rangle+\cdots+\langle h_n,h_n'\rangle,\;\;\;\;\|h\|=\sqrt{\langle h,h\rangle}.
\]
From now on, we will denote by $\S$ the unit sphere in $\CHd$ for this norm, namely
\[
\S=\{h\in\CHd:\|h\|=1\}.
\]

\subsection{The condition number}\label{sec:mu}
The (normalized) {\em condition number} at $(h,z)\in\CHd\times\Pn$ is defined as follows
\[
\mu(h,z)=\left\|h\right\|\,\left\|(Dh(z)\mid_{\,z^\perp})^{-1} \mbox{Diag}(\|z\|^{d_i-1}\sqrt{d_i})\right\|,
\]
or $\mu(h,z)=\infty$ if $Dh(\zeta)\mid_{\,z^\perp}$ is not invertible. Here, $\|h\|$ is the Bombieri-Weyl norm of $h$ and the second norm in the product is the operator norm . Note that assuming the system and a representative of a projective point are normalized, i.e., $h\in\S$ and $\|z\|=1$,
$$\mu(h,z) = \|(Dh(z)\mid_{\,z^\perp})^{-1} \mbox{Diag}(\sqrt{d_i})\|$$
is determined by the operator norm of the inverse of the Jacobian  $Dh(\zeta)$ restricted to the subspace orthogonal to $z$ multiplied with a diagonal matrix that makes formulas look nicer.

We also define the condition number of the system $h$ as
\begin{equation}\label{equ: condition of system}
\mu(h) = \sup_{z\in h^{-1}(0)} \mu(h,z).
\end{equation}
Note that with this definition $h$ is a regular polynomial system if and only if $\mu(h)<\infty$.

\subsection{A total degree homotopy}
For a real number $r>0$ (a common setting is $r=1$), let the start system be
\begin{equation}\label{eq:totdegreehom}
g=(X_1^{d_1}-r^{d_1}X_0^{d_1},\ldots,X_n^{d_n}-r^{d_n}X_0^{d_n}).
\end{equation}
This system is used in practice by various homotopy continuation software packages due to the simplicity of its evaluation and the simplicity of its solutions.

In \NAG\ a total-degree start system (with $r=1$) is created as follows:
\begin{Macaulay2}
\begin{verbatim}
i1 : R = CC[x,y,z];

i2 : f = {x^2+y^2+z^2, x*y}; -- target system

i3 : totalDegreeStartSystem f

        2    2   2    2
o3 = ({x  - z , y  - z },

      {(1, -1, 1), (-1, 1, 1), (1, 1, 1), (-1, -1, 1)})
\end{verbatim}
\end{Macaulay2}

\subsection{A conjecture by Shub and Smale}
The condition metric on the so-called {\em solution variety} $$V = \{(h,z)\in\S\times\Pn\ |\ h(z)=0 \}$$ is obtained by multiplying the metric inherited from the product $\S\times\Pn$ by the condition number $\mu$. The length of a homotopy path $t\rightarrow (h_t,\zeta_t)$ then equals
\begin{equation}\label{equ: length in condition metric}
\mathcal{C}_0(f,g,\zeta_0)=\int_0^T\mu(h_t,\zeta_t)\|\dot h_t,\dot \zeta_t\|\,dt.
\end{equation}

Consider the following initial pair:
\begin{equation}\label{equ: conjecture}
g=\begin{cases}\sqrt{d_1}X_0^{d_1-1}X_1\\ \vdots\\ \sqrt{d_n}X_0^{d_n-1}X_n\end{cases},\;\;\;\;e_0=\begin{pmatrix}1\\0\\ \vdots\\0\end{pmatrix}.
\end{equation}
The following is a slightly modified conjecture of Shub and Smale (the pair $(g,e_0)$ is a slightly modified version of the pair in the original conjecture) :
\begin{conjecture}\label{conj: Shub-Smale}
 \cite[Conjecture 2.4]{ShubSmale:Bezout-5-main-conjecture} For $g$ picked in the sphere $\S$ with the uniform probability distribution, the expected length of the homotopy path ${\rm E}(\mathcal{C}_0(f,g,e_0))$ is bounded by a polynomial in $N$.
\end{conjecture}

Proving polynomial average complexity for the linear homotopy with the initial pair above would give a {\em deterministic} algorithm settling Smale's 17th problem.
However, the tightest known bound on this complexity obtained so far is $N^{O(\log \log N)}$; see B\"urgisser and Cucker \cite{Burgisser-Cucker:near-Smale-17}. Most of the recent developments in this direction are inspired by \cite{Shub:Bezout-6-geodesics}.

In \NAG\ a {\em good initial pair} is created as follows:
\begin{Macaulay2}
\begin{verbatim}
i1 : R = CC[x,y,z];

i2 : f = {x^2+y^2+z^2, x*y}; -- target system

i3 : goodInitialPair f

o3 = ({1.41421x*z, 1.41421y*z}, {{0, 0, 1}})
\end{verbatim}
\end{Macaulay2}

Note that the start system $g$ of (\ref{equ: conjecture}) has $e_0$ as the only isolated solution with the smallest condition number possible, while the rest of zeros are described by $X_0=0$ and form a set isomorphic to $\P^{n-1}$.
In particular, system $g$ can't be used to compute all solutions of the target system $f$.

\subsection{Random linear homotopy}\label{subsec:random}

A randomized approach has been developed by Beltr\'an and Pardo \cite{Beltran-Pardo:probabilistic-Smale-17} (see also \cite{Beltran-Pardo:fast-linear-homotopy}): the basic idea is to construct an initial pair $(g,z_0)$ in a random fashion. First, pick a random start system $g$ in the sphere $\S$ with uniform distribution, then pick $z_0$ to be a random zero of $g$. While the first task is straightforward, a na\"ive approach to the second would be dependent on solving $g$. Nevertheless, there is a clever way to build a pair with the given properties pointed out in \cite{Beltran-Pardo:probabilistic-Smale-17} that depends only on the ability to pick random points in a sphere with the uniform probability distribution. The procedure is implemented in \NAG\ package along the description given in \cite{Beltran-Leykin:CHT}:
\begin{Macaulay2}
\begin{verbatim}
i1 : R = CC[x,y,z];

i2 : f = {x^2+y^2+z^2, x*y}; -- target system

i3 : randomInitialPair f

                               2
o3 : ({(- .11245 + .222955*ii)x  + (.28313 + .30427*ii)x*y + ... ,
                               2
       (.321684 + .0258085*ii)x  + (.200917 + .161032*ii)x*y + ... },

      {{.30545-.0745817*ii, .277793+.71328*ii, .550537-.11004*ii}})
\end{verbatim}
\end{Macaulay2}

Thus constructed random linear homotopy is shown to have average polynomial complexity in \cite{Beltran-Pardo:probabilistic-Smale-17} providing a uniform randomized algorithm that solves Smale's 17th problem.

\section{Certified linear homotopy}

Here we give a partial summary of the main constructions of \cite{Beltran-Leykin:CHT}.

\subsection{Linear homotopy}\label{sec:linear}
Given two systems $f$ and $g$ in the unit sphere $\S\subset\CHd$ we define a {\em linear homotopy} as a segment of the geodesic curve on $\S$ connecting $f$ and $g$ given by its arc--length parametrization:
\begin{equation}\label{equ: linear homotopy}
t\rightarrow h_t=g\cos(t)+\frac{f-\Real(\langle f,g\rangle) g}{\sqrt{1-\Real(\langle f,g\rangle)^2}}\sin(t),\;\;\;t\in\left[0,T\right],
\end{equation}
where
\[
T=\arccos\Real(\langle f,g\rangle)=\operatorname{distance}(f,g)\in[0,\pi].
\]

The procedure of certified tracking for a linear homotopy is presented by Algorithm~\ref{alg: certified homotopy}.

\begin{algorithm} \label{alg: certified homotopy} $z_*=TrackLinearHomotopy(f,g,z_0)$
\begin{algorithmic}[1]
\REQUIRE $f,g\in\S$; $z_0$ is an approximate zero of $g$.
\ENSURE $z_*$ is an approximate zero of $f$ associated to the end of the homotopy path starting at the zero of $g$ associated to $z_0$
and defined by the homotopy (\ref{equ: linear homotopy}).
\STATE $i \leftarrow 0$; $s_i=0$.
\WHILE {$s_i \neq T$}
\STATE Compute
\[
\dot g_i \leftarrow \dot h_s=-g\sin(s)+\frac{f-\Real(\langle f,g\rangle) g}{\sqrt{1-\Real(\langle f,g\rangle)^2}}\cos(s).
\]
at $s=s_i$.
\STATE $\varphi_i\leftarrow \chi_{i,1} \chi_{i,2}$ where
    \begin{eqnarray*}
    \chi_{i,1}&=&\left\|\binom{Dg_i(z_i)}{z_i^*}^{-1}\begin{pmatrix}\sqrt{d_1}& & & \\& \ddots & &\\& &\sqrt{d_n}\\&&&1\end{pmatrix}
\right\|\\
\chi_{i,2}&=&\sqrt{\|\dot g_i\|^2+\left\|\binom{Dg_i(z_i)}{z_i^*}^{-1}\binom{\dot g_i(z_i)}{0}\right\|^2}
    \end{eqnarray*}
\STATE Let $t_i$ be any number satisfying
\[
\frac{\LinearHomotopyConstant}{2d^{3/2}\varphi_i}\leq t_i\leq\frac{\LinearHomotopyConstant}{d^{3/2}\varphi_i}.
\]
\IF{$t_i>T-s_i$}
\STATE $t_i\leftarrow T-s_i$.
\ENDIF
\STATE $s_{i+1} \leftarrow s_i + t_i $; $g_{i+1} \leftarrow h_{s_{i+1}}.$
\STATE Perform a step of the projective Newton's method: $$z_{i+1} \leftarrow \|N_\P(g_{i+1})(z_{i})\|^{-1} N_\P(g_{i+1})(z_{i}).$$
\STATE $i \leftarrow i + 1$.
\ENDWHILE
\STATE $z_* \leftarrow z_T$.
\end{algorithmic}
\end{algorithm}

The correctness of the algorithm is shown in \cite{Beltran:proof-of-CHT}. The complexity of the algorithm is bounded by the following:
\[
|\, \text{projective Newton's method steps}\,| \leq \lceil 71d^{3/2}\mathcal{C}_0\rceil,
\]
where $d$ the maximal degree of polynomials in the system and $\mathcal{C}_0$ is defined by (\ref{equ: length in condition metric}).

The following is an example of usage of one of the main functions of \NAG\ with options that specify the computation to Algorithm~\ref{alg: certified homotopy}.
\begin{Macaulay2}
\begin{verbatim}
i1 : needsPackage "NumericalAlgebraicGeometry";

i2 : f = randomSd {2,2}; -- target system (random 2 quadrics)

i3 : (G,solution) = randomInitialPair f;

i4 : track(g, f, solution, Predictor=>Certified)

o4 = {{.050877+.571108*ii, -.74469-.020102*ii, -.154178-.304176*ii}}

i5 : (first oo).NumberOfSteps

o5 = 781
\end{verbatim}
\end{Macaulay2}
In this example a solution of a random system is found by a linear homotopy with a random initial pair in 781 steps.

\section{Experiments}

First we give outline several several experiments staged in \cite{Beltran-Leykin:CHT} and then present the results of an experimental search for optimal start systems.

\subsection{Average complexity comparison}\label{subsec: previous experiments}

We have obtained experimental data by running a linear homotopy connecting an initial pair $(g,z_0)$ to a random system in $\S\subset\CHd$ with $d_i=2$ for $i=1,\ldots,n$. Three kinds of initial pairs were examined:
\begin{itemize}
  \item {\bf good}: The initial pair (\ref{equ: conjecture}) conjectured to be ``good'' by Shub and Smale;
  \item {\bf total}: The start system $g$ of the total degree homotopy (\ref{eq:totdegreehom}) with $z_0 = (1,1,\ldots,1)$;
  \item {\bf random}: The random initial pair discussed in Subsection \ref{subsec:random}.
\end{itemize}

One can see the details of this experiment in~\cite{Beltran-Leykin:CHT}: to summarize, the conclusion is that the average complexity of {\bf good} is smaller than that of {\bf total}, which in turn is better than the complexity of {\bf random} (for all $n$, for which we ran the computation).

Here we report a stronger experimental {\bf negative result}: in experiments with the problems having at most 8 solutions we failed to find an initial pair that performed better than the ``good'' initial pair. (More precisely, a one-day long random search for each case with systems having at most 8 solutions fails to get an initial pair that performs better than~(\ref{equ: conjecture}) on average.)

This not only provides an experimental evidence to, but also prompts the following stronger version of the conjecture by Shub and Smale:
\begin{conjecture}\label{conj: stronger Shub-Smale}
  The initial pair $(g,e_0)$ as described in (\ref{equ: conjecture}) minimizes ${\rm E}(\mathcal{C}_0(f,g,z_0))$ over all initial pairs $(g,z_0)$ in the solution variety $V$.
\end{conjecture}

The weaker Conjecture~\ref{conj: Shub-Smale}, indeed, follows from Conjecture~\ref{conj: stronger Shub-Smale} in view of the existence result~\cite[Main Theorem]{ShubSmale:Bezout-5-main-conjecture}.

\subsection{Optimal start system}\label{subsec: new experiments}
Much like in the problem of finding one zero, in \cite{ShubSmale:Bezout-3} Shub and Smale use the condition number of a system as defined by (\ref{equ: condition of system}) as a guideline for picking a ``good'' start system. They proceed to show that in case of $n=1$ (one equation) the condition number of the total-degree start system (\ref{eq:totdegreehom}) grows exponentially as a function of the degree $d$ of the equation. In turn, the so-called {\em elliptic Fekete polynomials} give a family $F_d$ for which $\mu(F_d)=O(d)$. The construction of this family goes via Fekete points on the two-sphere, the problem of finding which relates to Smale's 7th problem.

To summarize, minimizing the condition number of a start system is a hard problem. The experiments below suggest that the complexity of the computation with a given start system, which ultimately is linked to the sum of the length of homotopy continuation paths in the condition metric
 $${\rm E}(\mathcal{C}_0(f,g)) = \sum_{z\in g^{-1}(0)}{\rm E}(\mathcal{C}_0(f,g,z)),$$
is not determined just by its condition number.

Our methodology is simple to describe (the choices of constants were made in such a way that the runtime of each of the computations reported in this article does not exceed one day):
\begin{enumerate}
  \item Among 1000000 random start systems $g$ picked on $\S$ with uniform probability distribution we select {\em five} with either the smallest $\mu(g)$ or the best estimated average complexity;
  \item Approximate the average number of steps the implementation of Algorithm~\ref{alg: certified homotopy} takes to compute all solutions by averaging over 10000 random target systems;
  \item Compare with the average performance of specially constructed start systems (total-degree, elliptic Fekete polynomials).
\end{enumerate}

For the construction of the total-degree systems with the minimal condition number we have created a script in \Mathematica\ that construct a (symbolic) expression for the condition number $\mu(g)$ depending on the parameter $r$ in (\ref{eq:totdegreehom}) and performs either symbolic or numerical optimization determining the optimal value for $r$.

\subsubsection{{\bf Case $n=1$}} 
\begin{figure}[h]
\begin{picture}(300,150)
   \put(3,0){\includegraphics[scale=1, keepaspectratio]{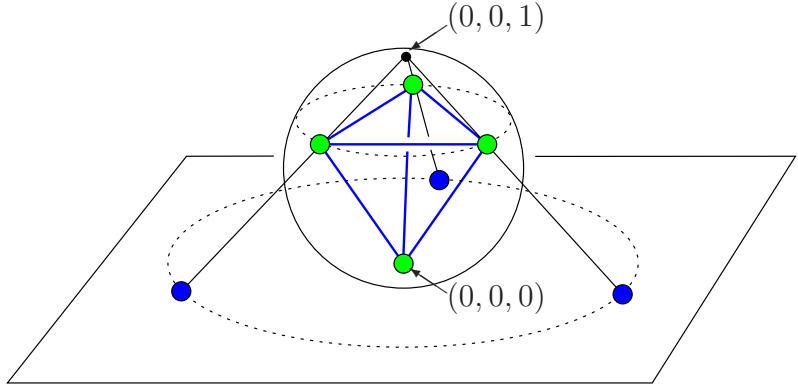}}
   \put(170,30){$(0,0,0)$}
   \put(170,35){\vector(-3,2){14}}
   \put(170,136){$(0,0,1)$}
   \put(170,136){\vector(-3,-2){14}}
\end{picture}
  \caption{Construction of an elliptic Fekete polynomial from four Fekete points.}
  \label{fig: Fekete4}
\end{figure}

Consider a single equation of degree four: one can look at
$$g_{total} = X_1^4-X_0^4,$$
as $r=1$ in (\ref{eq:totdegreehom}) is optimal for this case, and
$$g_{Fekete} = X_1(X_1^3-2\sqrt{2}X_0^3).$$
The construction of the latter can be carried out explicitly using the general recipe of \cite{ShubSmale:Bezout-3}. Four Fekete points are the vertices of the regular tetrahedron inscribed in the sphere $S^2$ of diameter $1$ and the center at $(0,0,\frac{1}{2})$. Placing one of the vertices in the origin (see Figure~\ref{fig: Fekete4}) the stereographic projections of the others onto the $xy$-plane lie on a circle of radius $\sqrt{2}$.

The results of two searches are reported below. Both use the methodology described in the beginning of the subsection, but differ in the way the systems are chosen in step (1).
\begin{itemize}
  \item The search presented in Figure~\ref{fig:one_equation_4_cond} is done according to the condition number $\mu(g)$;
  \item The search presented in Figure~\ref{fig:one_equation_4_step} uses the routine that estimates the average complexity.
\end{itemize}
The tables contain the detailed conditioning data as well as the estimated complexity (average number of steps in the certified homotopy continuation algorithm) for the specified start systems. The random systems $g_i$ in the tables are ordered according to $\mu(g_i)$.

\begin{figure}[h]
{\small
$$
\begin{array}{|r||c|c|c|c|c||c|c|}
  \hline
    & g_1 & g_2 & g_3 & g_4 & g_5 & g_{Fekete}& g_{total}\\ \hline
  \mu(g,z_1) & {1.23188}&
       {1.23094}&
       {1.23585}&
       {1.23975}&
       {1.24018}&
       {1.22475}&
       {1.41421}\\
   \mu(g,z_2) &     {1.23644}&
       {1.23788}&
       {1.23836}&
       {1.22948}&
       {1.20723}&
       {1.22475}&
       {1.41421}\\
    \mu(g,z_3) &    {1.21733}&
       {1.23569}&
       {1.19765}&
       {1.23153}&
       {1.23769}&
       {1.22475}&
       {1.41421}\\
     \mu(g,z_4) &   {1.23051}&
       {1.22567}&
       {1.23537}&
       {1.23855}&
       {1.23718}&
       {1.22474}&
       {1.41421}\\
       \hline
  \mu(g) &
       {1.23644}&
       {1.23788}&
       {1.23836}&
       {1.23975}&
       {1.24018}&
       {1.22475}&
       {1.41421}\\
  \hline
  \text{\#steps}  & {1138.32}&
       {1136.93}&
       {1137.81}&
       {1141.76}&
       {1130.81}&
       {1137.51}&
       {1177.79}\\
  \hline
\end{array}
$$
}
\caption{Search by condition number: $n=1$, $d=4$.}
\label{fig:one_equation_4_cond}
\end{figure}

\begin{figure}[h]
{\small
$$
\begin{array}{|r||c|c|c|c|c||c|c|}
  \hline
    & g_1 & g_2 & g_3 & g_4 & g_5 & g_{Fekete}& g_{total}\\ \hline
  \mu(g,z_1) & {1.28246}&
      {1.30382}&
      {1.33386}&
      {1.43745}&
      {1.44709}&
      {1.22475}&
      {1.41421}\\
  \mu(g,z_2) &  {1.25521}&
      {1.14504}&
      {1.3199}&
      {1.10411}&
      {1.34199}&
      {1.22475}&
      {1.41421}\\
   \mu(g,z_3) &    {1.23108}&
      {1.38061}&
      {1.35796}&
      {1.28299}&
      {1.25412}&
      {1.22475}&
      {1.41421}\\
    \mu(g,z_4) &   {1.17805}&
      {1.35972}&
      {1.14828}&
      {1.41302}&
      {1.25466}&
      {1.22474}&
      {1.41421}\\
      \hline
  \mu(g) &
      {1.28246}&
      {1.38061}&
      {1.35796}&
      {1.43745}&
      {1.44709}&
      {1.22475}&
      {1.41421}\\
  \hline
  \text{\#steps}  & {1138.47}&
      {1146.18}&
      {1152.21}&
      {1154.79}&
      {1155.42}&
      {1137.51}&
      {1177.79}\\
  \hline
\end{array}
$$
}
\caption{Search by average complexity: $n=1$, $d=4$.}
\label{fig:one_equation_4_step}
\end{figure}

We can see that both the elliptic Fekete polynomial and $g_i$ discovered by a random search perform significantly better than $g_{total}$.
A natural question arises that, to our knowledge, has not been answered yet: Is it true that the elliptic Fekete polynomials give optimal start equations in case $n=1$? The experimental data in Figures~\ref{fig:one_equation_4_cond}~and~\ref{fig:one_equation_4_step} provides no clear support to neither negative nor positive answer.

The phenomena that are observed in the experiments with one equation of degree 4 are amplified when the degree is increased. Figure~\ref{fig:one_equation_10_cond} contains the results for degree 10 which show that the random search (according to the condition number) gives much better performance both for the certified and heuristic homotopy tracking algorithm. (We used tighter than default tolerances in the heuristic runs in order to get the number of steps closer to that of the certified algorithm. See the corresponding \M2\ script for the exact settings used.)

\begin{figure}[h]
{\small
$$
\begin{array}{|r||c|c|c|c|c||c|c|}
  \hline
    & g_1 & g_2 & g_3 & g_4 & g_5 & g_{total}\\ \hline
  \mu(g) &
       {2.02342}&
      {2.02772}&
      {2.03071}&
      {2.03198}&
      {2.0452}&
      {7.15542}\\
  \hline
  \text{certified} &
      {11342}&
      {11395.2}&
      {11351}&
      {11393.6}&
      {11258.8}&
      {17737.2}\\
  \text{heuristic} &
  {2451.49}&
      {2475.89}&
      {2466.04}&
      {2441.97}&
      {2437.07}&
      {3376.32}\\
  \hline
\end{array}
$$
}
\caption{Search by condition number: $n=1$, $d=10$.}
\label{fig:one_equation_10_cond}
\end{figure}

\subsubsection{{\bf Case $n=2$}} One can show that in order to construct the best conditioned total-degree start system we should set $r \approx 0.746119$ in (\ref{eq:totdegreehom}). It is possible to express $r$ in radicals: a computation in Mathematica determines that $r$ has to be a root of $r^4(1 + 4 r^2)=1$.

Figure~\ref{fig:two_equations_2_cond} suggests that better start systems can be found via a random search.

\begin{figure}[h]
{\small
$$
\begin{array}{|r||c|c|c|c|c||c|c|}
  \hline
    & g_1 & g_2 & g_3 & g_4 & g_5 & g_{total}\\ \hline
  \mu(g) &
       {1.91848}&
      {1.92438}&
      {1.93129}&
      {1.93511}&
      {1.9386}&
      {2.23607}\\
  \hline
  \text{\#steps}  & {1259}&
      {1260.72}&
      {1265.15}&
      {1265.33}&
      {1283.8}&
      {1301.15}\\
  \hline
\end{array}
$$
}
\caption{Search by condition number: $n=2$, $d_1=d_2=2$.}
\label{fig:two_equations_2_cond}
\end{figure}

For this case already the search by average complexity becomes quite expensive with the current implementation of the certified tracking procedure and the computation of the size similar to that in case $n=1$ and $d=4$ produces systems much further from optimal than those in Figure~\ref{fig:two_equations_2_cond}.

As in the case $n=1$, larger examples exhibit larger gap in average complexity between the near-optimal systems found with the condition-number search and the total-degree start system. In case of two quartic equations, $n=2$ and $d_1=d_2=4$, the results are shown in Figure~\ref{fig:two_equations_4_cond}.

\begin{figure}[h]
{\small
$$
\begin{array}{|r||c|c|c|c|c||c|c|}
  \hline
    & g_1 & g_2 & g_3 & g_4 & g_5 & g_{total}\\ \hline
  \mu(g) & {4.18566}&
      {4.18963}&
      {4.204}&
      {4.32529}&
      {4.35056}&
      {4.91876}\\
  \hline
  \text{certified}  & {17003.6}&
      {17342.5}&
      {17143.8}&
      {17358.6}&
      {17209.6}&
      {20083.3}\\
  \text{heuristic}  & {8872.14}&
      {8782.57}&
      {8863.36}&
      {8944.19}&
      {8840.86}&
      {9599.85}\\
  \hline
\end{array}
$$
}
\caption{Search by condition number: $n=2$, $d_1=d_2=4$.}
\label{fig:two_equations_4_cond}
\end{figure}

We draw the following conclusions:
\begin{itemize}
  \item There is a way to find start systems that perform better than the total-degree start systems of the form (\ref{eq:totdegreehom}) by a random search.
  \item While bounding the condition number $\mu(g)$ of a system $g$ can be translated into good theoretical bounds (see \cite[Theorem~3.7]{Burgisser-Cucker:near-Smale-17}) on the average complexity, the latter does not depend solely on $\mu(g)$.
  \item If a system with the optimal $\mu(g)$ is explicitly known (as in the case $n=1$ and $d=4$), it is not clear whether it is optimal.
\end{itemize}

\section{Discussion}
In this article we have demonstrated the possibility of finding start systems that lead to better average complexity of finding all roots of a polynomial system than the systems of the total-degree family~(\ref{eq:totdegreehom}). Finding optimal start systems is a hard theoretical problem; however, good approximate solutions to this problem can be either constructed, e.g., by finding or approximating Fekete polynomials in case $n=1$,  or searched for experimentally. Once near-optimal start systems are available, these could be exploited in practice, e.g., in heuristic algorithms.

Here are some open questions, the answers to which we would like to know:
\begin{enumerate}
  \item If a start system $g$ has the optimal condition number $\mu(g)$ , does this imply that the average complexity of finding all roots (of a target system) with the corresponding linear homotopy is optimal? (This is similar to Conjecture~\ref{conj: stronger Shub-Smale} that deals with the case of just one root.)
  \item While the answer to the first question may turn out to be positive, one can show that $\mu(g_1)<\mu(g_2)$ does not force the same relation for the respective average complexities. However, what can be said if the relation between the condition numbers is uniform for all roots? I.e., if $\mu(g_1,z_1)<\mu(g_2,z_2)$ for all roots $z_1\in g_1^{-1}(0)$ and $z_2\in g_2^{-1}(0)$, is
      $$\sum_{z_1\in g_1^{-1}(0)} {\rm E}(\mathcal{C}_0(f,g_1,z_1))< \sum_{z_2\in g_2^{-1}(0)} {\rm E}(\mathcal{C}_0(f,g_2,z_2))\,?$$
  \item The near-optimal start systems we can find by a random search are not optimized for evaluation, which may be a bottleneck in the practical computation. Can one construct families of start systems that are, on one hand, easy to evaluate and, on the other hand, give better average complexity than the currently used families of total-degree start systems?
\end{enumerate}

The capabilities of the random search beyond cases with a small number of equations with a modest number of solutions are rather limited. Roughly speaking, this is due to the decreasing probability of picking a near-optimal system on the sphere $\S$ as the dimension of $\S$ grows. In view of this, a computationally inexpensive answer to question (3) above would be particularly valuable.

\bibliographystyle{abbrv}
\def\cprime{$'$}

\end{document}